\documentclass[12pt,draftcls,onecolumn]{IEEEtran}

\usepackage{bbold}
\usepackage{hyperref}
\usepackage{amsfonts}
\usepackage{amsmath,amssymb}
\usepackage{cite}
\usepackage{graphicx}
\usepackage{subcaption}
\usepackage{hyperref}
\usepackage{wrapfig}
\usepackage{mathrsfs}
\usepackage{color}

\usepackage{latexsym,url}
\usepackage{amsmath,amssymb,array,graphicx,verbatim}
\usepackage{setspace}
\usepackage{amsfonts}
\usepackage{times}

\newtheorem{theorem}{Theorem}

\newtheorem{corollary}[theorem]{Corollary}

\newtheorem{lemma}[theorem]{Lemma}

\newtheorem{proposition}[theorem]{Proposition}

\newcommand{\Real}{\mathbb{R}}
\newcommand{\Ex}{\mathbb{E}}
\newcommand{\Prob}{\mathbb{P}}
\newcommand{\base}{\mathbf{e}}
\newcommand\numberthis{\addtocounter{equation}{1}\tag{\theequation}}
\newcommand{\tr}{\ensuremath{{\scriptscriptstyle\mathsf{T}}}}

\newcommand{\R}{\mathbf{Cost}}
\newcommand{\V}{\mathbf{Var}}

\begin{document}

\title{Distributed Detection : Finite-time Analysis and Impact of Network Topology }

\author{Shahin Shahrampour$^\dagger$, Alexander Rakhlin$^\ddagger$ and Ali Jadbabaie$^\dagger$\footnote{

\noindent 
$\dagger$ Shahin Shahrampour and Ali Jadbabaie are with the Department of Electrical and Systems Engineering at the University of Pennsylvania, Philadelphia, PA 19104 USA. (e-mail: shahin@seas.upenn.edu; jadbabai@seas.upenn.edu).\\
$\ddagger$ Alexander Rakhlin is with the Department of Statistics at the University of Pennsylvania, Philadelphia, PA 19104 USA. (e-mail: rakhlin@wharton.upenn.edu).}} %

\maketitle

\begin{abstract}
This paper addresses the problem of distributed detection in multi-agent networks. Agents receive private signals about an {\it unknown} state of the world. The underlying state is globally identifiable, yet informative signals may be dispersed throughout the network. Using an optimization-based framework, we develop an iterative local strategy for updating individual beliefs. In contrast to the existing literature which focuses on asymptotic learning, we provide a {\it finite-time} analysis . Furthermore, we introduce a Kullback-Leibler cost to compare the efficiency of the algorithm  to its centralized counterpart.  Our bounds on the cost are expressed in  terms of network size, spectral gap, centrality of each agent and relative entropy of agents' signal structures. A key observation is  that distributing more informative signals to central agents results in a faster learning rate.  Furthermore, optimizing the weights, we can speed up learning by improving the spectral gap. We also quantify the effect of link failures on learning speed in symmetric networks. We finally provide numerical simulations which verify our theoretical results.
\end{abstract}

\section{Introduction}
Recent years have witnessed an intense interest on {\it distributed} detection, estimation, prediction and optimization \cite{tenney1981detection,tsitsiklis1993decentralized,borkar1982asymptotic,kar2012distributed,dekel2012optimal,nedic2009distributed,nedic2009distributed2}. Decentralizing the computation burden among agents has been widely regarded in networks ranging from sensor and robot to social and economic networks \cite{chamberland2003decentralized,bullo2009distributed,atanasov2014distributed,shahrampour2013online}.  
In this broad class of problems, agents in a network need to perform a global task for which they only have {\it partial} information. Therefore, they recursively exchange information with their neighbors, and the global dispersion of information in the network provides them with adequate data to accomplish the task. In the big picture, many of these schemes can also be embedded in the context of {\it consensus} protocols which have gained a growing popularity over the past three decades\cite{tsitsiklis1984problems,jadbabaie2003coordination,olfati2004consensus}. 

Earlier works on decentralized detection have considered scenarios where each agent sends its observations to a fusion center that decides over the true value of a parameter \cite{tenney1981detection,tsitsiklis1993decentralized,chamberland2003decentralized}. In these situations, the fusion center faces a classical hypothesis testing (centralized detection) problem after collecting the data from agents. Recently, another model of learning and detection has been proposed by Jadbabaie et al. \cite{jadbabaie2012non}. In this framework, the world is governed by a fixed true {\it state} or {\it hypothesis} that is aimed to be recovered by a network of agents. The state belongs to a {\it finite} set, and might represent a decision, an opinion, the price of a product or any quantity of interest. Each agent  observes a stream of {\it private} signals generated by a marginal of the global likelihood {\it conditioned} on the true state. However, the signals might {\it not} be informative enough for the agent to distinguish the underlying state of the world. Therefore, agents use {\it local} diffusion to compensate for their imperfect knowledge about the environment. In the literature, a host of schemes build on this model to describe distributed learning \cite{jadbabaie2012non,shahrampour2013exponentially,lalitha2014social,rahnama2010distributed}. Despite the wealth of results on the asymptotic behavior of these methods, the {\it finite-time} analysis remains elusive. In \cite{jadbabaie2012non}, a non-Bayesian update rule is proposed in the context of social networks. Each individual averages her Bayesian posterior belief with the opinion of her neighbors, and the beliefs tend to the truth under mild technical assumptions. Following up on the work of Duchi et al. \cite{duchi2012dual} on distributed dual averaging, an optimization-based algorithm is developed in \cite{shahrampour2013exponentially}. The authors demonstrate that the belief sequence generated according to their method is weakly consistent in undirected networks. Lalitha et al. \cite{lalitha2014social} introduce another strategy which puts exponential weights on a linear combination of Bayesian log-posteriors. The convergence conditions of their method are similar to those of \cite{jadbabaie2012non}. On the other hand, Rahnama Rad et al. \cite{rahnama2010distributed} present a distributed algorithm for {\it continuous} state space, and prove its convergence. In \cite{jadbabaie2012non,shahrampour2013exponentially,lalitha2014social}, the convergence occurs exponentially fast, and the {\it asymptotic} rate is characterized in terms of the {\it relative entropy} of individuals' signal structures and their {\it eigenvector centralities} (see \cite{jadbabaie2013information} for the rate analysis of \cite{jadbabaie2012non}). As an important consequence, the rate in \cite{shahrampour2013exponentially} only recovers the {\it empirical} average of relative entropies since the method is restricted to undirected networks.  

The asymptotic analysis presented in the above-discussed papers only describes the dominant factors that influence learning in the long run. In real world applications, however, the {\it decision} on the true state has to be made in a {\it finite} time. Therefore, it is crucial to study the finite-time variant of these schemes to gain insight into the interplay of {\it network parameters} which affect learning. To this end, we extend the work of Shahrampour et al. \cite{shahrampour2013exponentially} to directed networks where agents are {\it not} equally central. Moreover, we introduce the notion of {\it Kullback-Leibler (KL)} cost to measure the learning rate of an individual agent versus an {\it expert} who has all available information for learning. The KL decentralization cost simply compares the performance of distributed algorithm to its centralized counterpart. We derive an upper bound on the cost which proves the {\it spectral gap} of the network is substantial beside agents' centralities. It turns out that the upper bound scales inversely in the spectral gap, and logarithmically with the network size, number of states and time horizon. The rate also scales with the inverse of the  relative entropy of the conditional marginals. More specifically, the KL cost grows when signals do not provide enough evidence in favor of the true state versus some other state of the world.

Assuming that the network is realized with a {\it default} communication structure, each agent is endowed with a centrality. We establish that allocating more informative signals to more central agents can expedite learning. More interestingly, the importance of spectral gap opens new venues for {\it optimal} network design to facilitate agents' interactions. Each agent assigns different weights to its neighbors' information while communicating with them. We demonstrate how agents can modify these weights to achieve a faster learning rate. The key idea is to find the Markov chain with the best mixing behavior that is consistent with the network structure and agents' centralities. On the other hand, as a natural conjecture, we expect a more rapid learning rate in well-connected networks. We study the ramification of {\it link failures} in the network, and prove that in symmetric networks, less connectivity amounts to a sluggish learning process. We further apply our results on star, cycle and two-dimensional grid network. We observe that in each case the effect of spectral gap can be translated to the network {\it diameter}. Intuitively, a larger diameter makes the information propagation difficult around the network. Finally, we present numerical experiments which perfectly match our theoretical findings.

The rest of the paper is organized as follows: we describe the formal statement of the problem, and flesh out the distributed detection scheme in Section \ref{The Problem Description and Algorithm}. Section \ref{Finite-time Analysis of Beliefs and Cost Functions} is devoted to the finite-time analysis of the algorithm, whereas Section \ref{The Impact of Network Topology} elaborates on the impact of network characteristics on the convergence rate. We discuss briefly about applications of the model, and provide our numerical experiments in Section \ref{Numerical Experiment : Binary Signal Detection}. Section \ref{Conclusion} concludes. \\
\noindent
\textbf{Notation}: We adhere to the following notation in the exposition of our results:
\begin{center}
  \begin{tabular}{| c || l | }
    \hline
     $[n]$ &  The set $\{1,2,...,n\}$ for any integer $n$ \\ \hline
     $x^\tr$ & Transpose of the vector $x$ \\ \hline
     $x(k)$ & The $k$-th element of vector $x$ \\ \hline
     $x_{[k]}$ & The $k$-th largest element of vector $x$  \\ \hline
     $I_m$ &  Identity matrix of size $m$ \\ \hline
     $\Delta_m$ &  The $m$-dimensional probability simplex \\ \hline
     $\base_k $ & Delta distribution on $k$-th component \\ \hline
     $\left<\cdot, \cdot\right>$ &  Standard inner product operator \\ \hline
     $\|\cdot\|_p$ &  $p$-norm operator \\ \hline
     $\mathbf{1}$ &  Vector of all ones \\ \hline
     $\|\mu-\pi\|_{\text{TV}}$ &  Total variation distance between $\mu,\pi \in \Delta_m$ \\ \hline
     $D_{KL}(\mu \| \pi)$ &  KL-divergence of $\pi\in \Delta_m$ from $\mu \in \Delta_m$  \\ \hline
     $\lambda_i(W)$ &  The $i$-th largest eigenvalue of matrix $W$ \\ \hline
  \end{tabular}
\end{center}
For any $f\in \Real^m$ and $\mu \in \Delta_m$, we let $\Ex_\mu[\cdot]$ (respectively, $\V_\mu[\cdot]$) represent the expectation (respectively, variance) of $f$ under the measure $\mu$, i.e., we have 
\begin{align*}
\Ex_\mu[f]=\sum_{j=1}^m\mu(j)f(j) \ \ \ \ \ \ \ \ \ \ \V_\mu[f]=\sum_{j=1}^m\mu(j)\left(f(j)-\Ex_\mu[f]\right)^2.
\end{align*}

\section{The Problem Description and Algorithm}\label{The Problem Description and Algorithm}
In this section, we describe the observation and network model, and outline the centralized setting for the problem. Then, we provide a formal statement of the distributed setting, and characterize the distributing cost.

\subsection{Observation Model}
We consider an environment in which $\Theta=\{\theta_1,\theta_2,\ldots,\theta_m\}$ denotes a finite set of {\it states} of the world. We have a network of $n$ agents that seek the {\it unique}, true state of the world  $\theta_1\in \Theta$. At each time $t\in [T]$, the belief of agent $i$ is denoted by  $\mu_{i,t} \in  \Delta_m$, where $\Delta_m$ is a probability distribution over the set $\Theta$. In particular, $\mu_{i,0} \in \Delta_m$ denotes the prior belief of agent $i \in [n]$ about the states of the world assumed to be uniform with no loss of generality\footnote{The assumption of uniform prior only lets us avoid notational clutter. The analysis in the paper holds for any prior with full support.}. 

The learning model is given by a conditional likelihood function $\ell(\cdot|\theta_k)$ which is governed by a state of the world $\theta_k \in \Theta$. For each $i\in [n]$, let $\ell_i(\cdot|\theta_k)$ denote the $i$-th marginal of $\ell(\cdot|\theta_k)$, and we use the vector representation $\ell_i(\cdot|\theta)=[\ell_i(\cdot|\theta_1),...,\ell_i(\cdot|\theta_m)]^\tr$ to stack all states. At each time $t\in [T]$, the signal $s_t=(s_{1,t},s_{2,t},\ldots,s_{n,t})\in \mathcal{S}_1\times \dots \times \mathcal{S}_n$ is generated based on the true state $\theta_1$. Therefore, for each $i \in [n]$, the signal $s_{i,t}\in \mathcal{S}_i$ is a sample drawn according to the likelihood $\ell(\cdot|\theta_1)$ where $\mathcal{S}_i$ is the sample space. 

The signals are i.i.d. over time, and also the marginals are independent, i.e., $\ell(\cdot|\theta_k)=\Pi_{i=1}^n \ell_i(\cdot|\theta_k)$ for any $k \in [m]$. For the sake of convenience, we define $\psi_{i,t}\triangleq \log \ell_i(s_{i,t} |\theta)$ which is a sample corresponding to $\Psi_i\triangleq \log \ell_i(\cdot|\theta)$ for any $i\in [n]$.
\begin{description}
\item[{\bf A1.}] We assume that all log-marginals are uniformly bounded such that $\|\psi_{i,t}\|_\infty \leq B$ for any $s_{i,t}\in \mathcal{S}_i$, i.e., we have $|\log \ell_i(\cdot |\theta_k)| \leq B$ for any $i \in [n]$ and $k\in [m]$.
\end{description}
Assumption {\bf A1} is made for technical reasons, but such a bound can be found, for instance, when the signal space is discrete and provides a full support for distribution. Let us define $\bar{\Theta}_i$ as the set of states that are observationally equivalent to $\theta_1$ for agent $i\in [n]$; in other words, $\bar{\Theta}_i=\{\theta_k \in \Theta : \ell_i(s_i|\theta_k)= \ell_i(s_i| \theta_1) \ \ \forall s_i\in \mathcal{S}_i\} $ with probability one. As evident from the definition, any state $\theta_k \neq \theta_1$ in the set $ \bar{\Theta}_i$ is not distinguishable from the true state by observation of samples from the $i$-th marginal. Let $\bar{\Theta}=\cap_{i=1}^n\bar{\Theta}_i$ be the set of states that are observationally equivalent to $\theta_1$ from all agents perspective. 
\begin{description}
\item[{\bf A2.}] We assume that no state in the world is observationally equivalent to the true state from the standpoint of the network, i.e., the true state is globally identifiable, and we have $\bar{\Theta}=\{\theta_1\}$.  
\end{description}
Assumption {\bf A2} guarantees that the global likelihood provides sufficient information to make the true state uniquely identifiable. 

Let $\mathcal{F}_t$ be the smallest $\sigma$-field containing the information about all agents up to time $t$. Then, when the learning process continues for $T$ rounds, the probability triple $(\Omega,\mathcal{F},\Prob)$ is defined as follows: the sample space $\Omega=\otimes_{t=1}^T(\otimes_{i=1}^n\mathcal{S}_i)$, the $\sigma$-field $\mathcal{F}=\cup_{t=1}^{T}\mathcal{F}_t$, and the true probability measure $\Prob=\otimes_{t=1}^T \ell (\cdot| \theta_1)$. Finally, the operator $\Ex$ denotes the expectation with respect to $\Prob$.  

\subsection{Network Model}
The interaction between agents is captured  by a directed  graph $G=([n],E)$, where $[n]$ is the set of nodes corresponding to agents, and $E$ is the set of edges. Agents $i$ receives information from $j$ only if  the pair $(i,j) \in E$. We let $\mathcal{N}_i=\left\{j\in [n]: (i,j)\in E\right\}$ be the set of neighbors of agent $i$. Throughout the learning process agents truthfully report their information to their neighbors. We represent by $[W]_{ii}\geq 0$ the {\it self-reliance} of agent $i$, and by $[W]_{ij}>0$ the weight that agent $i$ assigns to information received from agent $j$ in its neighborhood. Then, the matrix $W$ is constructed such that $[W]_{ij}$ denotes the entry in its $i$-th row and $j$-th column. Therefore, $W$ has nonnegative entries, and $[W]_{ij}>0$ only if $(i,j)\in E$. For normalization purposes, we further assume that $W$ is stochastic; hence,
\begin{align*}
\sum_{j=1}^n[W]_{ij}=\sum_{j\in \mathcal{N}_i}[W]_{ij}=1.
\end{align*}
\begin{description}
\item[{\bf A3.}]We assume that the network is \textit{strongly connected}, i.e., there exists a directed path from any agent $i\in [n]$ to any agent $j\in [n]$. We further assume for simplicity that $W$ is diagonalizable\footnote{Note that the diagonalizability is not necessary, and it only forms a clean playground for technical analysis by avoiding Jordan blocks.}.
 \end{description}
The strong connectivity constraint in assumption {\bf A3} guarantees the information flow in the network. The assumption implies that $\lambda_1(W)=1$ is unique, and the other eigenvalues of $W$ are strictly less than one in magnitude \cite{rosenthal1995convergence}. Given the matrix of social interactions $W$, the eigenvector centrality is a non-negative vector $\pi$ such that for all $i\in [n]$,
\begin{align}\label{EigCent}
\pi(i)=\sum_{j=1}^n[W]_{ji}\pi(j).
\end{align}
for $\|\pi \|_1=1$. Then, $\pi(i)$ denoting the $i$-th element of $\pi$ is the eigenvector centrality of agent $i$. In the matrix form, the preceding relation takes the form $\pi^\tr W = \pi^\tr$, which means $\pi$ is the stationary distribution of $W$. Assumption {\bf A3} entails that the Markov chain $W$ is irreducible and aperiodic, and the unique stationary distribution $\pi$  has strictly positive components \cite{rosenthal1995convergence}.

\subsection{Centralized Detection}
To motivate the development of distributed scheme, we commence by introducing centralized detection\footnote{The method can be cast as special cases of  {\it Follow the Regularized Leader} \cite{abernethy2012interior} and {\it Mirror Descent} \cite{nemirovskii1983problem} algorithm.}. In this case, the scenario could be described as a two player repeated game between Nature and a {\it centralized} agent (expert) that has {\it global} information to learn the true state. More specifically, the expert observes the sequence of signals $\{s_t\}_{t=1}^T$ that are in turn revealed by Nature, and knows the entire network characteristics. At any round $t \in [T]$, the expert accumulates a {\it weighted average} of log-marginals, and forms the {\it belief} $\mu_t \in \Delta_m$ about the states, where $\Delta_m=\{\mu \in \Real^m \ | \ \mu \succeq 0, \ \sum_{k=1}^m\mu(k)=1\}$ denotes the $m$-dimensional probability simplex.  
Letting 
\begin{align}\label{psi definition} 
\psi_{t}\triangleq  \sum_{i=1}^n \pi(i)\psi_{i,t}=\sum_{i=1}^n \pi(i) \log \ell_i(s_{i,t}|\theta),
\end{align}
the sequence of interactions could be depicted in the form of the following algorithm:\\

{
\fbox{%
\centering
\begin{minipage}{5.9 in}
\textbf{Centralized Detection}\\
Input : A uniform prior belief $\mu_0$, a learning rate $\eta>0$.\\
Initialize : Let $\phi_0(k)=0$ for all $k\in [m]$.\\
At time $t=1,...,T$ : Observe the signal $s_t=(s_{1,t},s_{2,t},\ldots,s_{n,t})$, update the vector function $\phi_t$, and form the belief $\mu_t$ as follows,
\begin{align}
\phi_{t}=\phi_{t-1}+ \psi_{t}  \ \ \ \ \ \text{and} \ \ \ \ \ \mu_{t}&=\text{argmin}_{\mu \in \Delta_m} \left\{-\mu^\tr \phi_t + \frac{1}{\eta}D_{KL}(\mu \| \mu_{0}) \right\}  \label{CenBelief}.
\end{align}
\end{minipage}}}
\vspace{0.15 in}

At each time $t\in [T]$, the expert's goal is to maximize the expected log-marginals while sticking to the default belief $\mu_0$, i.e., minimizing the divergence. The trade-off between the two behavior is tuned with the {\it learning rate} $\eta$.

Let us note that according to Jensen's inequality for the concave function $\log(\cdot)$, we have for every $i\in [n]$
and $k\in [m]$ that
\begin{align*}
-D_{KL}\left(\ell_i(\cdot |\theta_1)\| \ell_i(\cdot  |\theta_k)\right)=\Ex\left[\log\frac{\ell_i(\cdot  |\theta_k)}{\ell_i(\cdot  |\theta_1)}\right] \leq \log\Ex\left[\frac{\ell_i(\cdot  |\theta_k)}{\ell_i(\cdot  |\theta_1)}\right]=0,
\end{align*}
where the inequality turns to equality if and only if $\ell_i(\cdot |\theta_1)= \ell_i(\cdot  |\theta_k)$, i.e., iff $\theta_k \in \bar{\Theta}_i$. Therefore, it holds that $\Ex[\log \ell_i(\cdot  |\theta_k)] \leq \Ex[\log \ell_i(\cdot  |\theta_1)]$, and recalling that the stationary distribution $\pi$ consists of positive elements, we have for any $k\neq 1$ that,
\begin{align*}
\Ex\left[\sum_{i=1}^n \pi(i) \Psi_i(k)\right]=\Ex\left[\sum_{i=1}^n \pi(i) \log \ell_i(\cdot  |\theta_k)\right] < \Ex\left[\sum_{i=1}^n \pi(i) \log \ell_i(\cdot  |\theta_1)\right]=\Ex\left[\sum_{i=1}^n \pi(i) \Psi_i(1)\right],
\end{align*}
where the strict inequality is due to uniqueness of the true state $\theta_1$, and the fact that $\bar{\Theta}=\cap_{i=1}^n\bar{\Theta}_i=\{\theta_1\}$ based on assumption {\bf A2}. In the sequel, without loss of generality, we assume the follwoing descending order, i.e.
 \begin{align}
\Ex\left[\sum_{i=1}^n \pi(i) \Psi_i(1)\right]>\Ex\left[\sum_{i=1}^n \pi(i) \Psi_i(2)\right]\geq \cdots \geq \Ex\left[\sum_{i=1}^n \pi(i) \Psi_i(m)\right], \label{orderassum}
\end{align}

\subsection{Distributed Detection}
We now extend the previous section to distributed setting modeled based on a network of agents. In the distributed scheme, each agent $i\in [n]$ only observes the stream of private signals $\{s_{i,t}\}_{t=1}^T$ generated based on the parametrized likelihood $\ell_i(\cdot|\theta_1)$. That is, agent $i\in [n]$ does not directly observe $s_{j,t}$ for any $j\neq i$. As a result, it gathers the local information by averaging the log-likelihoods in its neighborhood, and forms the belief $\mu_{i,t}\in \Delta_m$ at round $t\in [T]$ as follows:\\

{
\centering
\fbox{%
\begin{minipage}{6.4 in}
\textbf{Distributed Detection}\\
Input : A uniform prior belief $\mu_{i,0}$, a learning rate $\eta>0$.\\
Initialize : Let $\phi_{i,0}(k)=0$ for all $k\in [m]$ and $i\in [n]$.\\
At time $t\in [T]$ : Observe the signal $s_{i,t}$, update the function $\phi_{i,t}$, and form the belief $\mu_{i,t}$ as follows,
\begin{align}
\phi_{i,t}=\sum_{j\in \mathcal{N}_i} [W]_{ij}\phi_{j,t-1}+ \psi_{i,t}  \ \ \ \  \text{and} \ \ \ \ \mu_{i,t}&=\text{argmin}_{\mu \in \Delta_m} \left\{-\mu^\tr \phi_{i,t} + \frac{1}{\eta} D_{KL}(\mu \| \mu_{i,0}) \right\}  \label{DecBelief}.
\end{align}
\vskip 0.2in
\end{minipage}}

}
\vspace{0.15 in}

As outlined above, each agent updates its belief using purely local diffusion. We are interested in measuring the efficiency of the distributed algorithm via a metric comparing that to its centralized counterpart. At any round $t\in [T]$ , let us postulate that the cost which agent $i\in [n]$ needs to pay to have the same opinion as the expert is $D_{KL}(\mu_{i,t}\| \mu_t)$; then, the total {\it decentralization cost} that the agent incurs after $T$ rounds is as follows
\begin{align}\label{Regret}
\R_{i,T}\triangleq \sum_{t=1}^T D_{KL}(\mu_{i,t}\| \mu_t)= \sum_{t=1}^T \Ex_{\mu_{i,t}}\left[\log \frac{\mu_{i,t}}{ \mu_t}\right].
\end{align}
The function quantifies the difference between the agent that observes private signals $\{s_{i,t}\}_{t=1}^T$ and an expert that has $\{s_t\}_{t=1}^T$ and $\pi$ available. Note importantly that $\R_{i,T}$ is a random quantity since the expectation is {\it not} taken with respect to randomness of signals.

We conclude this section with the following lemma which reiterates that both algorithms are reminiscent of the well-known {\it Exponential Weights} algorithm.

\begin{lemma}\label{Beliefs}
The update rules \eqref{CenBelief} and \eqref{DecBelief} have the explicit form solutions, 
\begin{align*}
\mu_t(k)=\frac{\exp\{\eta \phi_t(k)\}}{\left<\mathbf{1},\exp\{\eta \phi_t\}\right>} \ \ \ \ \ \text{and} \ \ \ \ \ \mu_{i,t}(k)=\frac{\exp\{\eta \phi_{i,t}(k)\}}{\left<\mathbf{1},\exp\{\eta \phi_{i,t}\}\right>},
\end{align*} 
respectively, for any $i \in [n]$ and $k\in[m]$. Moreover,
\begin{align*}
\phi_{i,t}= \sum_{\tau=1}^{t}\sum_{j=1}^n\left[W^{t-\tau}\right]_{ij}\psi_{j,\tau}.
\end{align*}
\end{lemma}
We will now state the main results of the paper with underlying intuition behind them. The proofs are sometimes omitted and provided later in the appendix.

\section{Finite-time Analysis of Beliefs and Cost Functions}\label{Finite-time Analysis of Beliefs and Cost Functions}
In this section, we investigate the convergence of agents' beliefs to the true state in the network. Agents exchange information over time, and reach consensus about the true state. The connectivity of the network plays an important role in the learning as $W^t \rightarrow \mathbf{1}\pi^\tr/n$ as $t\rightarrow \infty$. To examine the learning rate, we need to have knowledge about the mixture behavior of Markov chain $W$. The following lemma sheds light on the mixture rate, and we invoke it later for technical analysis. 
\begin{lemma}\label{mixture}
Given strong connectivity of the network (assumption {\bf A3}), the stochastic matrix $W$ satisfies
\begin{align*}
\sum_{\tau=1}^{t}\sum_{j=1}^n\left|\left[W^{t-\tau}\right]_{ij}-\pi(j)\right| \leq \frac{4\log n}{1-\lambda_{\max}(W)},
\end{align*}
for any $i \in [n]$, where $\lambda_{\max}(W)\triangleq \max \left\{\left|\lambda_{n}(W)\right|,\left|\lambda_{2}(W)\right|\right\}$
\end{lemma}
We now establish that agents have arbitrarily close opinions in a connected network. Furthermore, the convergence rate is governed by cardinality of state space and network characteristics. 
\begin{lemma}\label{Distributed Beliefs}
Let the sequence of beliefs $\{\mu_{i,t}\}_{t=1}^T$ for each agent $i\in [n]$ be generated by the Distributed Detection algorithm with the learning rate $\eta$. Given bounded log-marginals (assumption {\bf A1}), global identifiability of the true state (assumption {\bf A2}), and strong connectivity of the network (assumption {\bf A3}), for each individual agent $i\in [n]$ it holds that
\begin{align*}
\frac{1}{\eta}\log \|\mu_{i,t}-\base_1\|_{\text{TV}} &\leq - \mathcal{I}(\theta_1,\theta_2) t  + \sqrt{2B^2t\log \frac{m}{\delta}} + \frac{8B\log n}{1-\lambda_{\max}(W)}+\frac{\log m}{\eta},
\end{align*}
with probability at least $1-\delta$, where for $k\geq 2$
\begin{align*}
\mathcal{I}(\theta_1,\theta_k)\triangleq \sum_{i=1}^n\pi(i)D_{KL}(\ell_i(\cdot|\theta_1)\|\ell_i(\cdot|\theta_k)).
\end{align*}
\end{lemma}
Lemma \ref{Distributed Beliefs} verifies that the belief $\mu_{i,t}$ of each agent $i\in [n]$ is {\it strongly} consistent, i.e., it converges almost surely to a delta distribution on the true state. The claim follows immediately by letting $\delta=1/t^2$ and applying Borel-Cantelli lemma. However, we are interested in the interplay of parameters in finite-time and in particular the behavior of decentralization cost function in \eqref{Regret}. Let us now proceed to the next lemma to derive a variance-type bound on the cost.

\begin{lemma}\label{DecRegretDecay}
The decentralization cost function \eqref{Regret} associated to the Distributed Detection algorithm with the learning rate $\eta$ satisfies
\begin{align*}
\R_{i,T} \leq  2\eta^2\sum_{t=1}^T \V_{\mu_t}\left[q_{i,t}\right], 
\end{align*}
so long as $\eta\|q_{i,t}\|_{\infty}\leq 1/4$ at each round, where $q_{i,t}\triangleq \phi_{i,t}-\phi_{t}$.
\end{lemma}
The bound in Lemma \ref{DecRegretDecay} is evocative of numerous regret bounds developed for the well-known problem of prediction with expert advice corresponding to the centralized detection in an adversarial setting \cite{cesa2006prediction,rakhlin2013online}. However, such bounds are in terms of second moment, rather than variance which is a smaller quantity (see e.g. the bound in Lemma 3 of \cite{rakhlin2013online} derived in terms of local norms). The following theorem illuminates how the variance bound comes in handy by concentrating the measure around the true distribution.
\begin{theorem}\label{DecRegretRate}
Let the sequence of beliefs $\{\mu_{i,t}\}_{t=1}^T$ for each agent $i\in [n]$ be generated by the Distributed Detection algorithm with the choice of learning rate $\eta=\frac{1-\lambda_{\max}(W)}{16B\log n}$. Given bounded log-marginals (assumption {\bf A1}), global identifiability of the true state (assumption {\bf A2}), and strong connectivity of the network (assumption {\bf A3}), we have
\begin{align*}
 \R_{i,T} \leq \max \left\{  \frac{8B^2}{\mathcal{I}^2(\theta_1,\theta_2)}\log\left[\frac{mT}{\delta}\right] , \frac{4B\log n}{\mathcal{I}(\theta_1,\theta_2)}\frac{\log \left[mT\right]}{1-\lambda_{\max}(W)} \right\} +1,
\end{align*}
with probability at least $1-\delta$.
\end{theorem} 
\begin{IEEEproof} We recall that $q_{i,t}$ in the statement of Lemma \ref{DecRegretDecay} satisfies
\begin{align*}
\|q_{i,t}\|_\infty &= \left\| \sum_{\tau=1}^{t}\sum_{j=1}^n\left(\left[W^{t-\tau}\right]_{ij}-\pi(j)\right)\psi_{j,t} \right\|_\infty
&\leq B\sum_{\tau=1}^{t}\sum_{j=1}^n\left|\left[W^{t-\tau}\right]_{ij}-\pi(j)\right|
&\leq \frac{4B\log n}{1-\lambda_{\max}(W)},
\end{align*}
due to Lemma \ref{mixture} and assumption {\bf A1}. Therefore, the choice of $\eta=\frac{1-\lambda_{\max}(W)}{16B\log n}$ guarantees that $q_{i,t}$ satisfies $\eta\|q_{i,t}\|_{\infty}\leq 1/4$ for all $t\in [T]$. We now explicitly calculate the variance of $q_{i,t}$ under the measure $\mu_t$. Then, we apply H$\ddot{\text{o}}$lder's inequality for primal-dual norm pairs to bound it as,
\begin{align*}
\V_{\mu_t}\left[q_{i,t}\right] &= \sum_{k=1}^m \mu_t(k)\left(q_{i,t}(k)-\Ex_{\mu_t}\left[q_{i,t}\right]\right)^2\\
&= \sum_{k=1}^m \mu_t(k)\left(\left<q_{i,t},\base_k\right>-\left<q_{i,t},\mu_t\right>\right)^2\\
& \leq \left<q_{i,t},\base_1-\mu_t\right>^2 + \sum_{k=2}^m \mu_t(k)\left<q_{i,t},\base_k-\mu_t\right>^2\\
&\leq \big\|q_{i,t}\big\|^2_\infty \big\|\base_1-\mu_t \big\|^2_1 + \sum_{k=2}^m \mu_t(k)\big\|q_{i,t}\big\|^2_\infty \big\|\base_k-\mu_t \big\|^2_1\\
&\leq \big\|q_{i,t}\big\|^2_\infty \big\|\base_1-\mu_t \big\|^2_1 + 4\big\|q_{i,t}\big\|^2_\infty \sum_{k=2}^m \mu_t(k)\\
&= 4\big\|q_{i,t}\big\|^2_\infty \big\|\base_1-\mu_t \big\|^2_{\text{TV}} + 4\big\|q_{i,t}\big\|^2_\infty \big\|\base_1-\mu_t \big\|_{\text{TV}} ,
\end{align*}
where in the last line we used the fact that $\|\base_k-\mu_t \|_1=2\|\base_k-\mu_t \|_{\text{TV}}\leq 2$ for any $k \in [m]$. Taking into account the condition $\eta \|q_{i,t}\|_{\infty}\leq 1/4$, we obtain
\begin{align}
2\eta^2 \V_{\mu_t}\left[q_{i,t}\right] \leq \frac{1}{2 }\left(\big\|\base_1-\mu_t \big\|^2_{\text{TV}}+\big\|\base_1-\mu_t \big\|_{\text{TV}}\right) \leq \big\|\base_1-\mu_t \big\|_{\text{TV}}. \label{Prod1}
\end{align}
Following exactly the same steps in the proof of Lemma \ref{Distributed Beliefs}, it can be verified that for any $t\in [T]$, the centralized algorithm yields
\begin{align*}
\frac{1}{\eta}\log \|\mu_{t}-\base_1\|_{\text{TV}} &\leq -\mathcal{I}(\theta_1,\theta_2) t + \sqrt{32B^2t\log \frac{m}{\delta}} +\frac{\log m}{\eta},
\end{align*}
with probability at least $1-\delta$. To have the identity above work for every $t \in [T]$ with probability at least $1-\delta$, we need to take a union bound over all $t \in [T]$, which changes the parameter $\delta$ to $\delta/T$ in the right hand side of the preceding relation. Let us avoid notational clutter, by defining $a\triangleq \mathcal{I}(\theta_1,\theta_2)$ and $b\triangleq (32B^2\log\left[mT/\delta\right])^{1/2}$, respectively. Then, in view of the identity above, with probability at least $1-\delta$ we can bound \eqref{Prod1} as follows,
\begin{align*}
2\eta^2 \V_{\mu_{t}}\left[q_{i,t}\right] &\leq  m \exp \left\{- a\eta  t+ b\eta \sqrt{t} \right\}\\
&\leq  m \exp \left\{-\frac{a}{2} \eta t \right\} \ \ \ \ \text{for} \ \ \ \ t \geq t_1\triangleq \left(\frac{2b}{a}\right)^2\\
&\leq  \frac{1}{T}   \ \ \ \ \text{for} \ \ \ \ t \geq t_2 \triangleq \frac{2}{a\eta}\log \left[mT\right]. 
\end{align*} 
Let $t_0=\max \{t_1,t_2\}$ and consider the relation in above as well as the condition $\eta \|q_{i,t}\|_{\infty}\leq 1/4$ to observe
\begin{align*}
2\sum_{t=1}^T \eta^2 \V_{\mu_{t}}\left[q_{i,t}\right]  &= 2\sum_{t=1}^{t_0} \eta^2 \V_{\mu_{t}}\left[q_{i,t}\right] + 2\sum_{t=t_0+1}^{T} \eta^2 \V_{\mu_{t}}\left[q_{i,t}\right]\\
&\leq 2\sum_{t=1}^{t_0} \Ex_{\mu_t}[\eta^2 q^2_{i,t}] + \sum_{t=t_0+1}^{T} \frac{1}{T}\\
&\leq  2 \sum_{t=1}^{t_0} \frac{1}{16} + 1 = \frac{t_0}{8} +1,
\end{align*}
with probability at least $1-\delta$. Plugging the bound above into Lemma \ref{DecRegretDecay} completes the proof.
\end{IEEEproof}
Regarding Theorem \ref{DecRegretRate} the following comments are in order: the rate is related to the inverse of $\mathcal{I}(\theta_1,\theta_2)$ which is a weighted average of KL-divergence of observations under $\theta_2$ (the second best alternative) from observations under $\theta_1$ (the true state). Also, from the definition of $\mathcal{I}(\theta_1,\theta_2)$ in Lemma \ref{Distributed Beliefs}, the weights turn out to be agents' centralities. Intuitively, when signals hardly reveal the difference between the best two candidates for the true state, agents must make more effort to distinguish the two. In turn, this results in suffering a larger cost caused by slower learning. 
The decentralization cost always scales logarithmically with the number of states $m$. Now define 
\begin{align}\label{spectral gap}
\gamma(W)\triangleq 1-\lambda_{\max}(W), 
\end{align}
as the {\it spectral gap} of the network. Then, Theorem \ref{DecRegretRate} suggests that for large networks, the cost scales inversely in the spectral gap, and logarithmically with the network size $n$. Finally, the detection cost with respect to time horizon is $\mathcal{O}(\log T)$ which is sub-linear. Therefore, the average cost (per iteration cost) asymptotically tends to zero. Moreover, such dependence is quite natural as even an expert incurs a $\mathcal{O}(\log T)$ regret to detect the true state\cite{cesa2006prediction}.

\section{The Impact of Network Topology}\label{The Impact of Network Topology}
The results of previous section verify that network characteristics govern the learning process. We now discuss the role of agents' centralities and the network spectral gap. 
\subsection{Effect of Agent Centrality}
To examine centrality, let us return to the definition of $\mathcal{I}(\theta_1,\theta_2)$ in Lemma \ref{Distributed Beliefs}, and imagine that the network is {\it collaborative} in the sense that the network designer wants to expedite learning. Then, to have the best information dispersion, the marginal which collects the most evidence in favor of $\theta_1$ against $\theta_2$ should be allocated to the most central agent. By the same token, in an {\it adversarial} network where Nature aims to delay the learning process, such marginal should be assigned to the least central agent. To sum up, let us put forth the concept of network {\it regularity} as defined in \cite{jadbabaie2013information} in the context of social learning. Recalling the definition of eigenvector centrality \eqref{EigCent}, we say a network $G$ is more regular than $G'$ if $\pi'$ majorizes $\pi$, i.e., if for all $j\in [n]$
\begin{align}\label{Regularity}
\sum_{i=1}^j \pi_{[i]}  \leq \sum_{i=1}^j \pi'_{[i]},   
\end{align}
where $\pi_{[i]}$ denotes the $i$-th largest element of $\pi$. Letting 
\begin{align*}
u\triangleq \left[D_{KL}(\ell_1(\cdot|\theta_1)\|\ell_1(\cdot|\theta_2)),\ldots,D_{KL}(\ell_n(\cdot|\theta_1)\|\ell_n(\cdot|\theta_2))\right]^\tr,
\end{align*}
it is a straightforward consequence of Lemma 1 proved in \cite{jadbabaie2013information} that 
\begin{align*}
\sum_{i=1}^n \pi_{[i]}u_{[i]}  \leq \sum_{i=1}^n \pi'_{[i]}u_{[i]},
\end{align*} 
when $\pi'$ majorizes $\pi$. Therefore, spreading more informative signals among central agents speeds up the learning procedure. 

\subsection{Optimizing the Spectral Gap}

We now turn our attention to the spectral gap of network \eqref{spectral gap}. Suppose that agents are given a default communication matrix $W$ which determines their neighborhood and centrality. The problem is to find the optimal spectral gap assuming that the neighborhood and centrality of each agent are fixed. The key idea is to change the mixing behavior of the Markov chain $W$. It is well-known, for instance, that we could do so using lazy random walks \cite{levin2009markov} which replaces $W$ with $\frac{1}{2}(W+I_n)$. To generalize the idea, let us define a modified communication matrix
\begin{align}\label{MixS}
W' \triangleq \alpha W + (1-\alpha)I_n  \ \ \ \alpha \in [0,1],
\end{align}
which has the same eigenstructure as $W$. Then, the eigenvalues of $W'$ are weighted averages of those of $W$ with one. From standpoint of network designing, one can exploit the freedom in choosing $\alpha$ to optimize the spectral gap. 
\begin{proposition}\label{OPT}
The optimal spectral gap of the modified communication matrix $W'$ \eqref{MixS} is as follows,
\begin{align*}
\gamma^*= \frac{2-2\lambda_2(W)}{2-\lambda_n(W)-\lambda_2(W)}  \ \ \ \text{for} \ \ \ \alpha^*=\frac{2}{2-\lambda_n(W)-\lambda_2(W)},
\end{align*}
when $\lambda_n(W)+\lambda_2(W)<0$
\end{proposition}
\begin{IEEEproof}
To optimize the spectral gap, we need to minimize the second largest eigenvalue of $W'$ in magnitude, that is, to solve the min-max problem
\begin{align}\label{minmax}
\min_{\alpha\in [0,1]} \lambda_{\max}(W')=\min_{\alpha\in [0,1]} \max\left\{|\alpha\lambda_2(W)+1-\alpha|,|\alpha\lambda_n(W)+1-\alpha|\right\}.
\end{align}
Drawing the plots of $|\alpha\lambda_2(W)+1-\alpha|$ and $|\alpha\lambda_n(W)+1-\alpha|$ in terms of $\alpha$ verifies that the minimum occurs at the intersection of the lines 
\begin{align*}
\alpha\lambda_2(W)+1-\alpha=-\alpha\lambda_n(W)+\alpha-1,
\end{align*}
yielding $\alpha^*=\frac{2}{2-\lambda_n(W)-\lambda_2(W)}$. Plugging $\alpha^*$ into the min-max problem \eqref{minmax}, we calculate the optimal value $\lambda_{\max}^*$ as
\begin{align*}
\lambda_{\max}^*=\frac{\lambda_2(W)-\lambda_n(W)}{2-\lambda_n(W)-\lambda_2(W)},
\end{align*}
and since $\gamma^*=1-\lambda_{\max}^*$ the proof follows immediately.
\end{IEEEproof}

\subsection{Sensitivity to Link Failure}
It is intuitive that in a network with more links, agents are offered more opportunities for communication. Adding  links provides more avenues for spreading information, and improves the learning quality. We study this phenomenon for symmetric networks where a pair of agents assign similar weights to each other, i.e., $W^\tr=W$. In particular, we explore the connection of spectral gap with the link failure. In this regard, let us introduce the following positive semi-definite matrix
\begin{align}\label{PertMatrix}
\Delta W(i,j)\triangleq (\mathbb{e}_i-\mathbb{e}_j)(\mathbb{e}_i-\mathbb{e}_j)^\tr,
\end{align}
where $\mathbb{e}_i$ is the $i$-th unit vector in the standard basis of $\Real^n$. Then, for $i,j \in [n]$ the matrix 
\begin{align}\label{NewComm}
\bar{W}(i,j)\triangleq W+ [W]_{ij} \Delta W(i,j),
\end{align}
corresponds to a new communication matrix that removes edges $(i,j)$ and $(j,i)$ from the network, and adds $[W]_{ij}=[W]_{ji}$ to the self-reliance of agent $i$ and agent $j$.
\begin{proposition}\label{LF}
Consider the communication matrix $\bar{W}(i,j)$ in \eqref{NewComm}. Then, for any $i,j\in [n]$ the following identity holds
\begin{align*}
\lambda_{\max}\left(W\right) \leq \lambda_{\max}\left(\bar{W}(i,j)\right),
\end{align*}
so long as $W$ is positive semi-definite.
\end{proposition}
\begin{IEEEproof}
We recall that $\Delta W(i,j)$ in \eqref{PertMatrix} is positive semi-definite with $\lambda_n\left(\Delta W(i,j)\right)=0$. Applying Weyl's eigenvalue inequality on \eqref{NewComm}, we obtain for any $k\in [n]$
\begin{align*}
\lambda_k\left(W\right) \leq \lambda_k\left(\bar{W}(i,j)\right),
\end{align*}
which holds in particular for $k=2$. On the other hand, the matrix $W$ is positive semi-definite, so we have that $\lambda_{\max}\left(W\right)=\lambda_2\left(W\right)$. Combining with the fact that $\bar{W}(i,j)$ is symmetric and positive semi-definite, the proof is completed.
\end{IEEEproof}
The proposition immediately implies that removing a link reduces the spectral gap. In this case, in view of the bound in Theorem \ref{DecRegretRate}, the decentralization cost has more latitude to vary. Therefore, to keep the costs small, agents tend to maintain their connections. Let us take note of the delicate point that monotone increase in the upper bound does {\it not} necessarily imply a monotone increase in the cost; however, one can {\it roughly} expect such behavior. We elaborate on this issue in the numerical experiments. Finally, notice that the positive semi-definiteness constraint on $W$ is not strong, since it can be easily satisfied by replacing a lazy random walk $\frac{1}{2}(W+I_n)$ with $W$.

\subsection{Star, Cycle and Grid Networks}
\begin{figure}[ht!]
\centering
  \begin{tabular}{@{}ccc@{}}
    \includegraphics[scale=1.65]{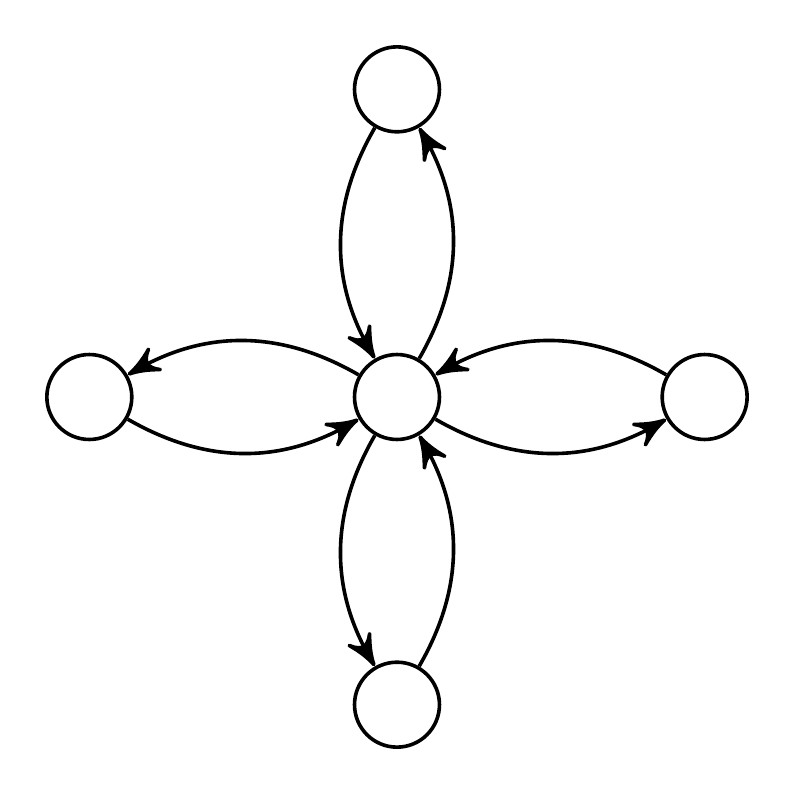} &
    \includegraphics[scale=1.1]{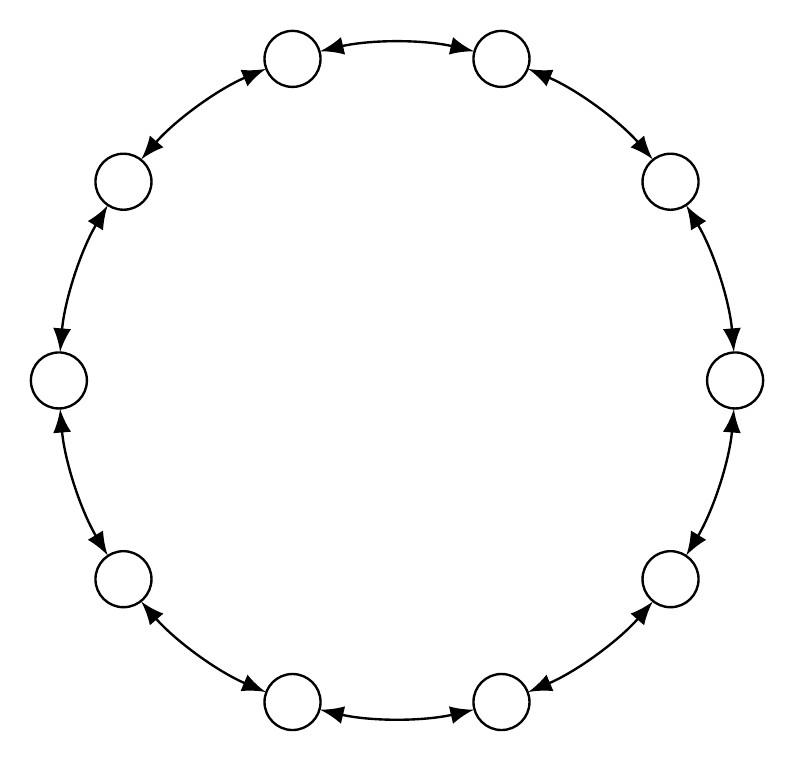} &
    \includegraphics[scale=1.1]{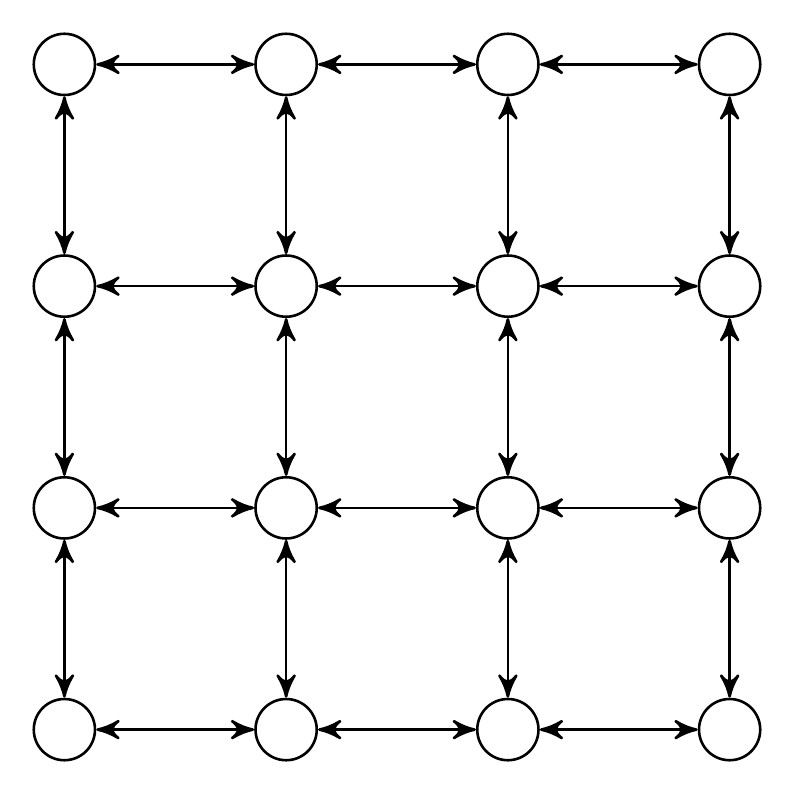} 
  \end{tabular}
\caption{Illustration of networks :  star, cycle and grid networks with $n$ agents. For each network, each individual agent possesses a self-reliance of $\omega \in (0,1)$.}
\label{Example of Nets}
\end{figure}

We now examine the spectral gap impact for some interesting networks (Fig. \ref{Example of Nets}), and derive explicit bounds for decentralization cost. As one of the famous examples in computer networks, we start with the star network. Regardless of the network size, existence of one central agent always preserves the network diameter, and therefore, we expect a benign scaling with network size. On the other side of the spectrum lies the cycle network where the diameter grows linearly with the network size. We should, hence, observe how the poor communication in cycle network affects the learning rate. Finally, as a possible model for sensor networks, we study the grid network where the network size scales quadratically with the diameter.

\begin{corollary}\label{examples}
Under conditions of Theorem \ref{DecRegretRate} and the choice of learning rate $\eta=\frac{\gamma(\cdot)}{16B\log n}$, for $n$ large enough we have the following bounds on the decentralization cost:
\begin{itemize}
\item[{\bf (a)}] For the star network in Fig. \ref{Example of Nets}
\begin{align*}
\R_{i,T} \leq \mathcal{O}\left(\frac{\log\left[nmT\right]}{\min\left\{1-\omega , 1-|2\omega-1|\right\}}\right).
\end{align*}
\item[{\bf (b)}] For the cycle network in Fig. \ref{Example of Nets}
\begin{align*}
\R_{i,T} \leq   \mathcal{O}\left(\frac{\log\left[nmT\right]}{\min\left\{1-|2\omega-1| , 2(1-\omega)\sin^2\frac{\pi}{n}\right\}}\right).
\end{align*}
\item[{\bf (c)}] For the grid network in Fig. \ref{Example of Nets}
\begin{align*}
\R_{i,T} \leq  \mathcal{O}\left(\frac{\log\left[nmT\right]}{\min\left\{1-|2\omega-1| , 2(1-\omega)\sin^2\frac{\pi}{\sqrt{n}}\right\}}\right).
\end{align*}
\end{itemize}
\end{corollary}
\begin{IEEEproof}
The spectrum of the Laplacian of star and cycle graphs are well-known \cite{chung1997spectral}. We have the eigenvalue set corresponding to communication matrix of star and cycle graphs as
\begin{align*}
\bigg\{1,\omega,\ldots,\omega,2\omega-1\bigg\} \ \ \ \text{and} \ \ \ \bigg\{\omega+(1-\omega)\cos\frac{2\pi i}{n}\bigg\}_{i=0}^{n-1},
\end{align*} 
respectively. Therefore, the proof of {\bf (a)} and {\bf (b)} follows immediately. The grid graph is the Cartesian product of two rings of size $\sqrt{n}$ (due to wraparounds at the edges), and hence, its eigenvalues are derived by summing the eigenvalues of two $\sqrt{n}$-rings\cite{chung1997spectral}. 
Therefore, the eigenvalue set takes the form 
\begin{align*}
\bigg\{\omega+(1-\omega)\cos\frac{\pi (i+j)}{\sqrt{n}}\cos\frac{\pi (i-j)}{\sqrt{n}}\bigg\}_{i,j=0}^{\sqrt{n}-1},
\end{align*}
and the proof of {\bf (c)} is completed. 
\end{IEEEproof}
Let us use the notation $\tilde{\mathcal{O}}(\cdot)$ to hide the poly log factors. Then, the bounds derived in Corollary \ref{examples} indicate that the algorithm requires $\tilde{\mathcal{O}}(1)$ iterations to achieve a near optimal log-distance from the true state in the star network. However, the rate deteriorates to $\tilde{\mathcal{O}}(n^2)$$($respectively, $\tilde{\mathcal{O}}(n))$ in the cycle (respectively, grid) network. In all cases, the rate is proportional to the diameter of the network which is a natural indicator of information dissemination quality.

\section{Numerical Experiment : Binary Signal Detection}\label{Numerical Experiment : Binary Signal Detection}
We now discuss distributed detection of signals transmitted through noisy channels. We first particularize the model to binary signals, and then present our simulation results in that context. 
\subsection{Signal Detection in Communication Channels}
\begin{figure}[t!]
\centering
\includegraphics[trim = 30mm 53mm 30mm 45mm, clip, scale=0.5]{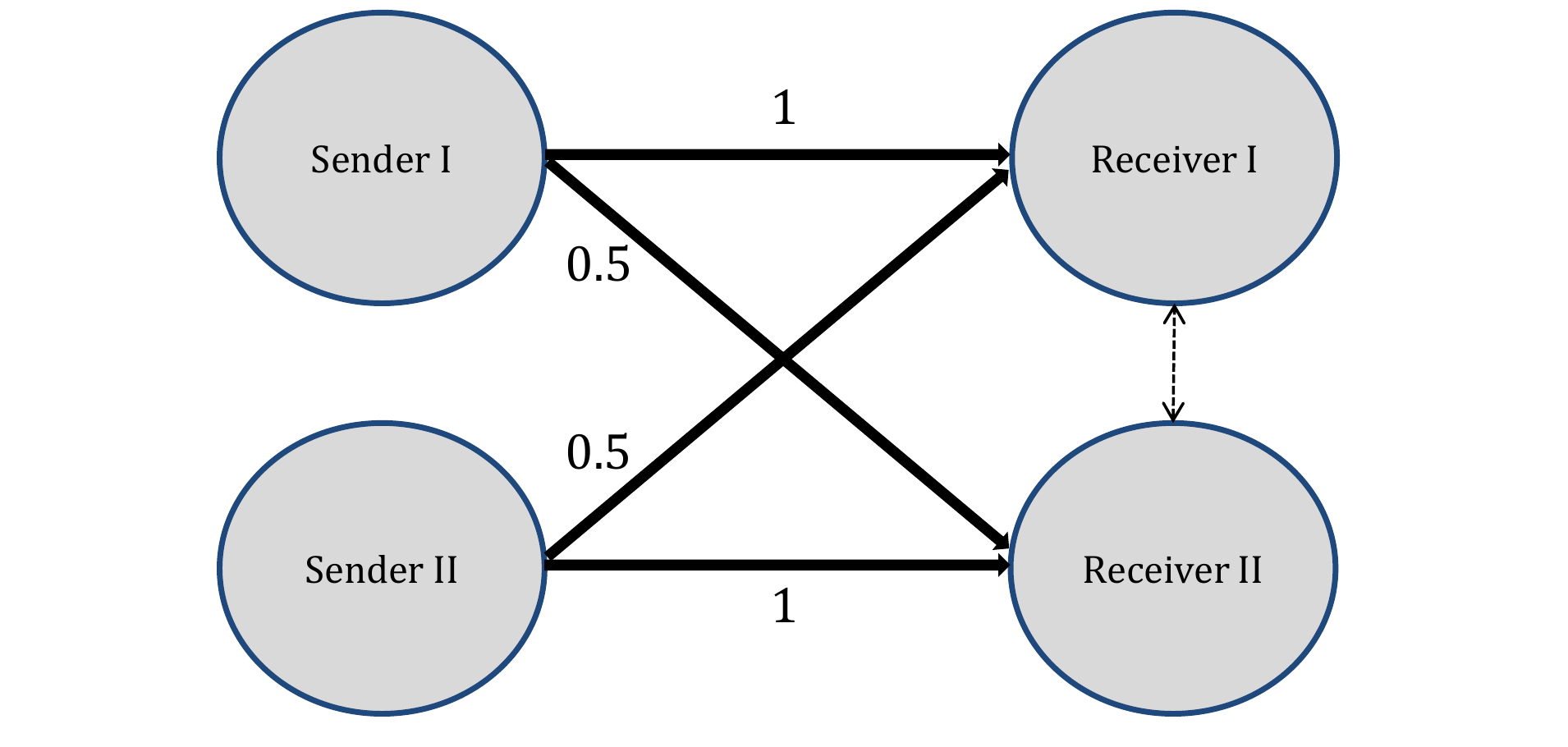}
\caption{A communication channel which transmits digital data. Each receiver cannot distinguish the message based on its own signals, so it communicates with the other receiver to identify the message.}
\label{CChannel}
\end{figure}

In information theory, data transmission can be modeled via a sender, a receiver and a channel. The channel is used to convey information from one end to another. In general, a faulty communication is possible, and it might be caused by channel noises, and imperfect {\it modulation} or {\it demodulation} (see e.g. \cite{cover2012elements,poor1994introduction}). In what follows, we exemplify this point, and employ distributed detection to resolve it.

Suppose a 2-digit binary number is to be transmitted over a communication channel as depicted in Fig. \ref{CChannel}. Sender {\bf I} and sender {\bf II} broadcast $T$ copies of the first and second digit, respectively. Receiver {\bf I} (agent {\bf I}) can recognize the first digit accurately\footnote{
To have the assumption {\bf A1} satisfied, we can think of accurate transmission as $1-\varepsilon$ probability of success for some small $\varepsilon>0$.} while the second digit is distorted with probability 1/2. On the other hand, receiver {\bf II} (agent {\bf II}) collects the exact value of the second digit at the terminal, and observes a misrepresented first digit with probability 1/2. In this example, the state space is $\Theta=\{\theta_1=00,\theta_2=01,\theta_3=10,\theta_4=11\}$, and let the true state be $\theta_1=00$. We can see that none of the receivers can solely establish a reliable communication with senders as each of them has difficulty inferring one digit. More formally, it is straightforward to calculate that
\begin{align*}
\ell_1(s_1| 00 )=\ell_1(s_1| 01) \ \ \text{$\forall s_1 \in \{0,1\}^2$} \ \ \ \ \ \ \text{and} \ \ \ \ \ \ \ell_2(s_2| 00 )=\ell_2(s_2| 10) \ \ \text{$\forall s_2 \in \{0,1\}^2$},
\end{align*}  
which simply means $\bar{\Theta}_1=\{\theta_1,\theta_2\}$ and $\bar{\Theta}_2=\{\theta_1,\theta_3\}$. However, the global identifiability of the true state holds as we have $\bar{\Theta}=\bar{\Theta}_1 \cap \bar{\Theta}_2 = \{\theta_1\}$. Therefore,  according to Lemma \ref{Distributed Beliefs}, exchanging information with each other, receivers are able to decipher the message transmitted by senders.

\subsection{Convergence of Beliefs}
\begin{figure}[t!]
\centering
\includegraphics[trim = 20mm 84mm 20mm 84mm, clip, scale=0.5]{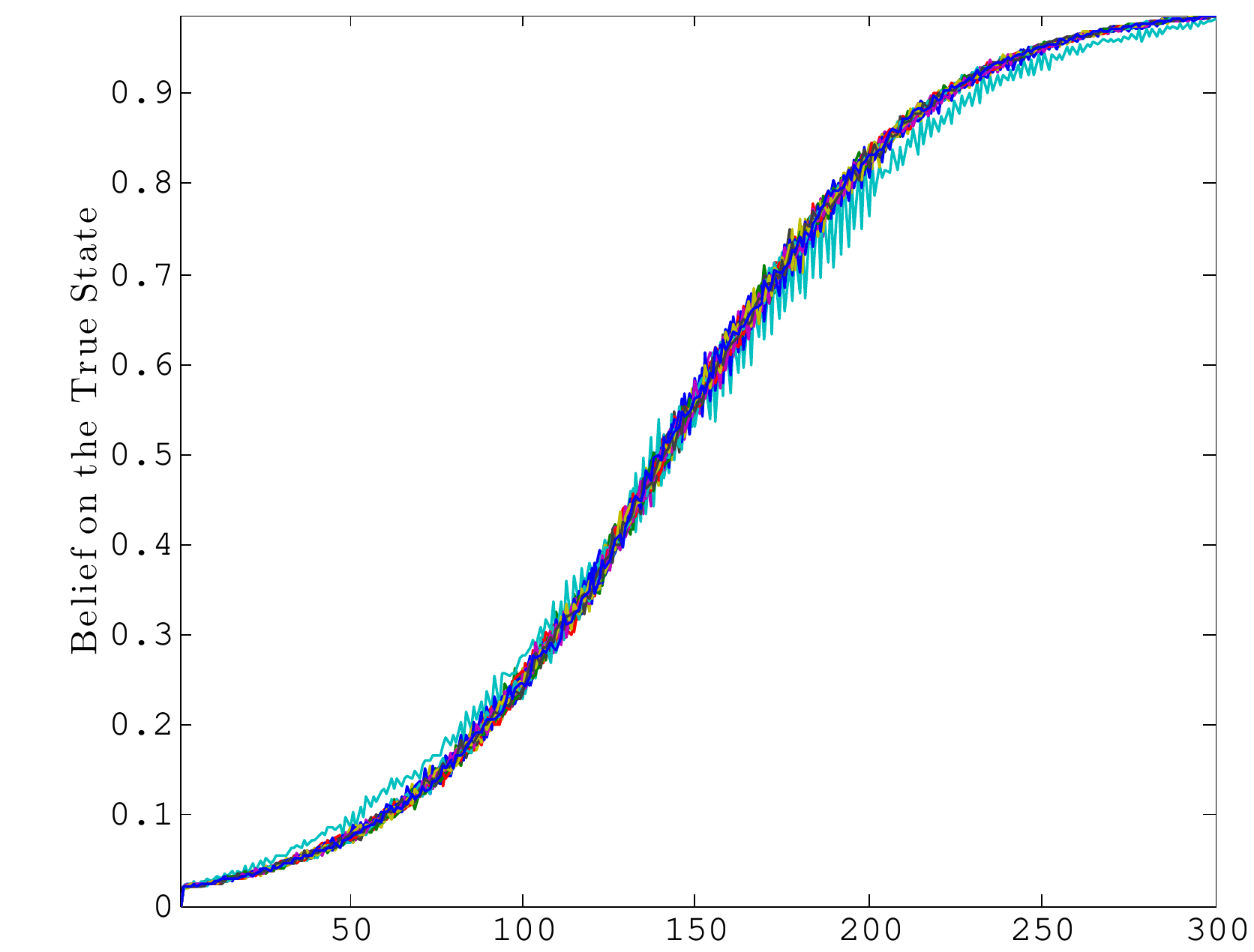}
\caption{The belief evolution for all 50 agents in the network. The global identifiability of the true state and strong connectivity of the network result in learning.}
\label{Convergence}
\end{figure}

For purpose of simulation, we generate a strongly connected network of $n=50$ agents with a default communication matrix $W$. Assume that there exist $m=51$ states in the world and agents are to discover the true state $\theta_1$. At time $t\in [T]$, a signal $s_{i,t}\in \{0,1\}$ is generated based on the true state such that $\ell_i(\cdot|\theta_1)=\ell_i(\cdot|\theta_{i+1})$. In other words, for agent $i\in [n]$, we have $\bar{\Theta}_i=\{\theta_1,\theta_{i+1}\}$ and $\theta_{i+1}$ is observationally equivalent to the true state. Therefore, each agent $i\in [n]$ fails to distinguish $\theta_1$ from $\theta_{i+1}$ once relying on the private signals. However, since we have $\bar{\Theta}=\cap_{i=1}^n\bar{\Theta}_i=\{\theta_1\}$, the true state is globally identifiable. Consequently, in view of Lemma \ref{Distributed Beliefs}, we expect that all agents reach a consensus on the true state (Fig. \ref{Convergence}), and learn the truth exponentially fast. 

\subsection{Optimizing the Spectral Gap}
\begin{figure}[ht!]
\centering
\includegraphics[trim = 8mm 74mm 9mm 82mm, clip, scale=0.6]{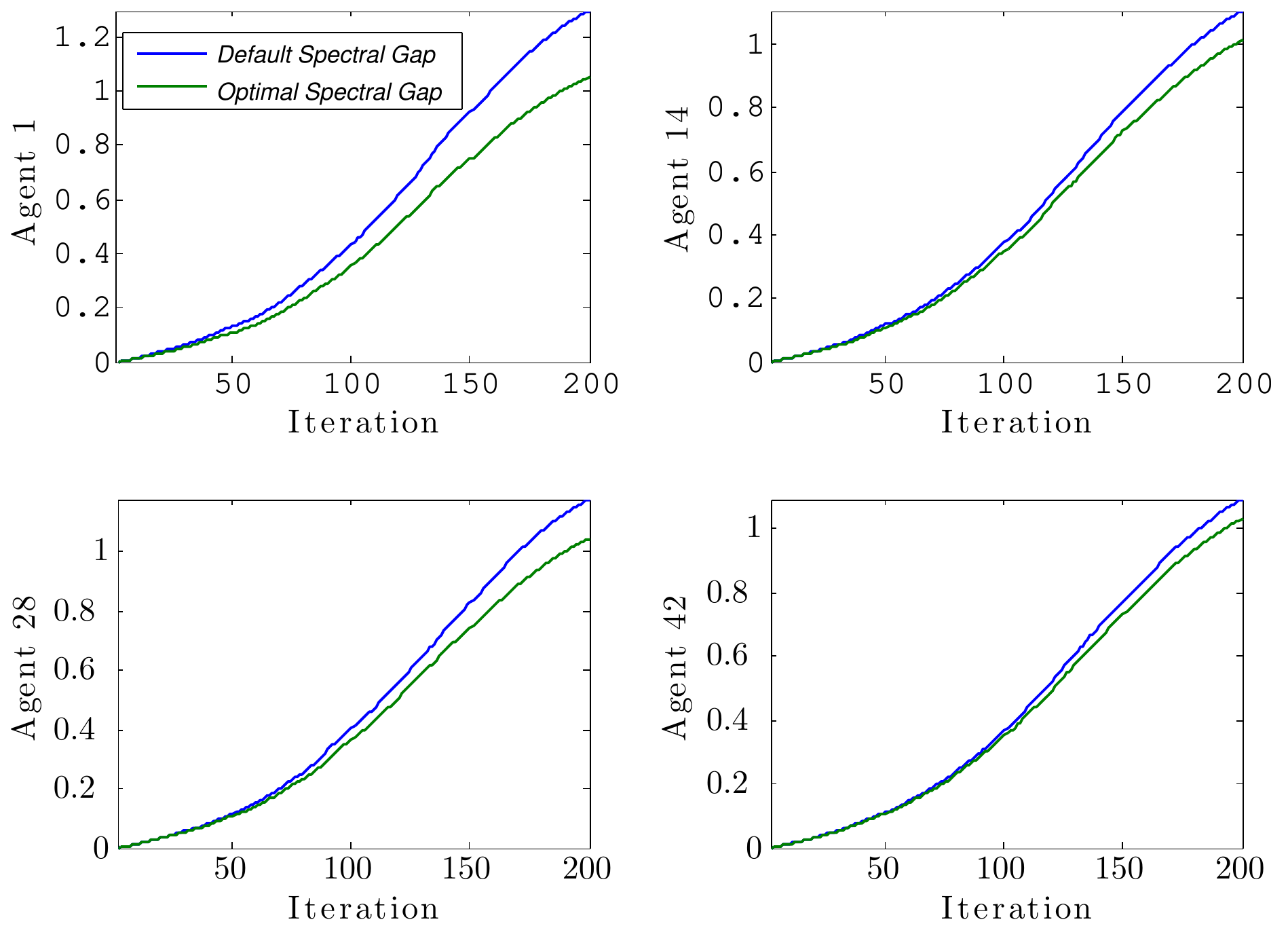}
\caption{The plot of decentralization cost versus time horizon for agents 1, 14, 28 and 42 in the network. The cost in the network with the optimal spectral gap (green) is always less than the network with default weights (blue).}
\label{DecentralizationCost}
\end{figure}

We now turn to optimizing the spectral gap to speed up learning. We proved in Proposition \ref{OPT} that every default communication matrix can be adjusted to a matrix $W'$ which has the optimal spectral gap when centralities are fixed. Setting the parameter $\alpha$ in \eqref{MixS} equal to $\alpha^*$ derived in Proposition \ref{OPT}, we obtain the optimal network. The dependence of decentralization cost to the spectral gap was theoretically proved in Theorem \ref{DecRegretRate}. Applying the results of  Proposition \ref{OPT} verifies that in the optimal network, agents suffer a lower decentralization cost comparing to the default network (Fig. \ref{DecentralizationCost}).

\subsection{Sensitivity to Link Failure}
\begin{figure}[ht!]
\centering
\includegraphics[trim = 0mm 79mm 3mm 84mm, clip, scale=0.6]{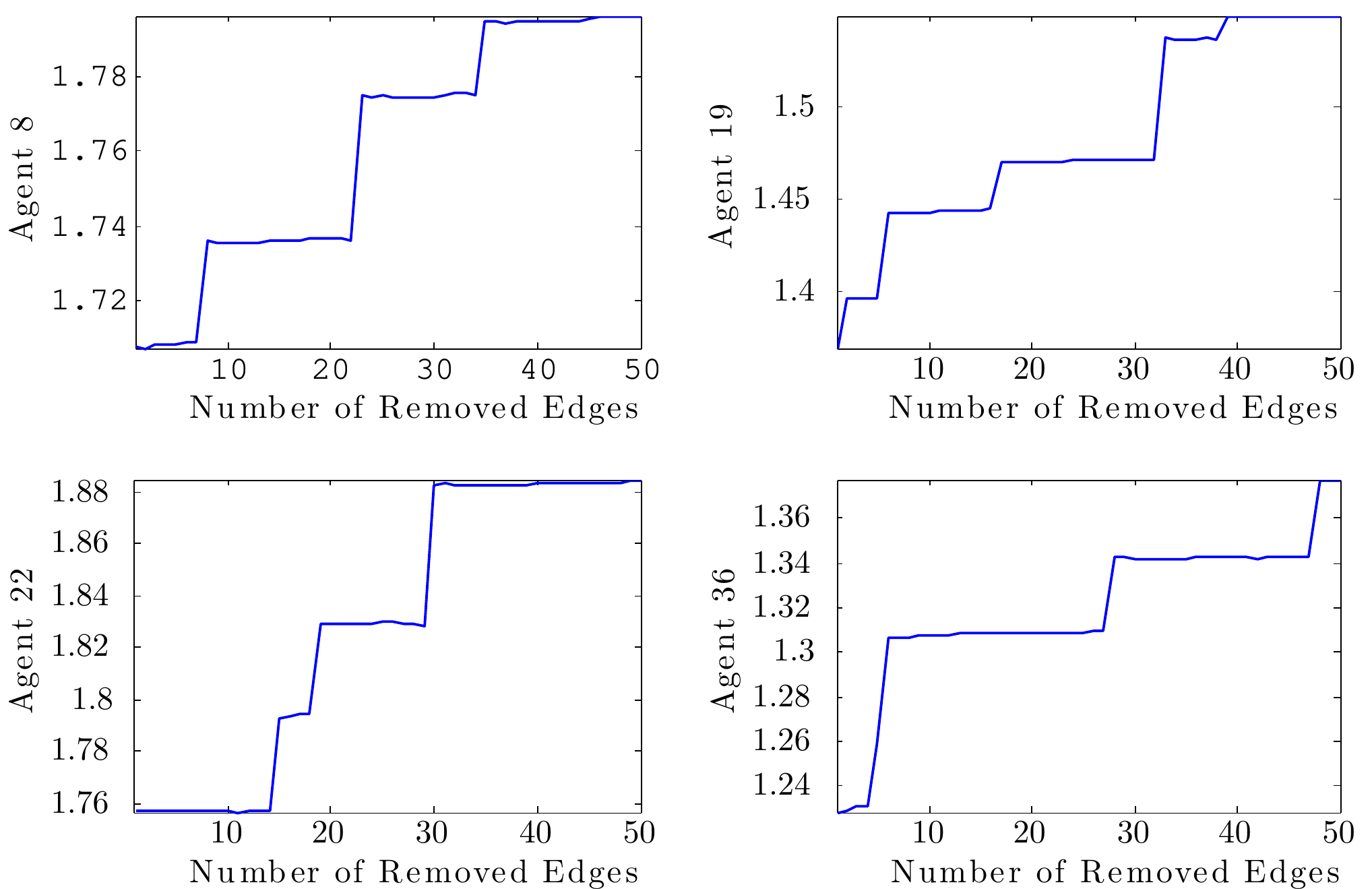}
\caption{The decentralization cost at round $T=300$ for agents 8, 19, 22 and 36 in the network. Removing the links causes poor communication among agents and increase the decentralization cost.}
\label{RemovingEdges}
\end{figure}

Let us symmetrize the network in the previous section such that $[W]_{ij}=[W]_{ji}$. In this case every agent is equally central, and we have $\pi=\mathbf{1}/n$. To study the impact of link failure, we sequentially select a random pair of agents in the network, and remove their connection. Each time that a link is discarded, we compute the decentralization cost in the new network at iteration $T=300$, and continue the process until 50 bi-directional edges are eliminated from the network. In view of Proposition \ref{LF}, we expect a monotone decrease in the spectral gap which amounts to a larger decentralization cost. We plot the cost for four agents in the network, and observe that the behavior is almost (not quite) monotonic (Fig. \ref{RemovingEdges}). The monotone dependence of the upper bound to the spectral gap (Theorem \ref{DecRegretRate}) does not necessarily guarantee a monotone relationship between cost and the spectral gap. Therefore, we can only roughly expect such behavior.

\section{Conclusion}\label{Conclusion}
We considered a distributed detection model where a network of agents aim to learn the underlying state of the world. The private signals do not provide enough information for agents about the true state. Hence, agents engage in a local communication to compensate for their imperfect knowledge. Each agent iteratively forms a belief about the state space using the collected data in its neighborhood. We analyzed the learning procedure for a {\it finite} time horizon. To study the efficiency of our algorithm versus its centralized counterpart, we brought forward the idea of KL cost. It turned out that network size, spectral gap, centrality of each agent and relative entropy of agents' signal structures are the key parameters that affect distributed detection. We established that allocating more informative signals to central agents as well as optimizing the spectral gap can speed up learning. We also proved that the learning rate deteriorates in the case of link failures, which can be seen as a side effect of poor communication. Finally, we would like to address a few issues in future works. In this paper, we discussed a communication model in which agents exchange information at every round. In some networks, all-time communication is potentially costly or unnecessary. Alternatively, agents can only contact each other when their signals are not informative enough about the true state. As another direction, we can consider scenarios where the signal distributions are not stationary. This generalizes the model to dynamic parameters where we can investigate detection robustness in changing environments.  

\section*{Appendix : Proofs}

\textbf{\emph{Proof of Lemma \ref{Beliefs}}}. The proof is elementary, and it is only given to keep the paper self-contained. We write the Lagrangian associated to the update \eqref{CenBelief} as,
\begin{align*}
L(\mu,\lambda)= -\mu^\tr \phi_t + \frac{1}{\eta}\left<\mu,\log\frac{\mu}{\mu_0}\right>+\lambda \mu^\tr\mathbf{1}-\lambda,
\end{align*}  
where we left the positivity constraint implicit. Differentiating above with respect to $\mu$ and $\lambda$, and setting the derivatives equal to zero, we get
\begin{align*}
\mu_{t}(k)=\mu_{0}(k)\exp\left\{\eta \phi_t(k)-\lambda-1\right\}\ \ \ \ \ \text{and} \ \ \ \ \ \mu_{t}^\tr\mathbf{1}=1,
\end{align*}
respectively, for any $k\in [m]$. Combining the equations above and noting that $\mu_{0}$ is uniform, we have
\begin{align*}
\frac{1}{m}\exp\{-\lambda-1\}\sum_{k=1}^m\exp\{ \eta \phi_t(k)\}=1,
\end{align*}
which allows us to solve for $\lambda$ and calculate the optimal solution $\mu_{t}$ as follows,
\begin{align*}
\mu_{t}(k)=\frac{\exp\left\{\eta \phi_t(k)\right\}}{\sum_{k=1}^m\exp\left\{\eta \phi_t(k)\right\}}.
\end{align*}
The proof for $\mu_{i,t}$ follows precisely in the same fashion. To calculate $\phi_{i,t}$, notice that in view of the first update in \eqref{DecBelief} we have
\begin{align*}
\left[ \begin{array}{cccc}
\phi_{1,t} \\
\phi_{2,t} \\
\vdots \\
\phi_{n,t}  \end{array} \right]=\left(W\otimes I_m\right)\left[ \begin{array}{cccc}
\phi_{1,t-1} \\
\phi_{2,t-1} \\
\vdots \\
\phi_{n,t-1}  \end{array} \right]+\left[ \begin{array}{cccc}
\psi_{1,t} \\
\psi_{2,t} \\
\vdots \\
\psi_{n,t}  \end{array} \right],
\end{align*}
where $\otimes$ denotes the Kronecker product. The equation above represents a discrete-time linear system. Given the fact that $\phi_{i,0}(k)=0$ for all $k\in [m]$ and $i \in [n]$, the closed-form solution of the system takes the form
\begin{align*}
\left[ \begin{array}{cccc}
\phi_{1,t} \\
\phi_{2,t} \\
\vdots \\
\phi_{n,t}  \end{array} \right]=\sum_{\tau=1}^t \left(W\otimes I_n\right)^{t-\tau}\left[ \begin{array}{cccc}
\psi_{1,\tau} \\
\psi_{2,\tau} \\
\vdots \\
\psi_{n,\tau}  \end{array} \right]=\sum_{\tau=1}^t \left(W^{t-\tau}\otimes I_n\right)\left[ \begin{array}{cccc}
\psi_{1,\tau} \\
\psi_{2,\tau} \\
\vdots \\
\psi_{n,\tau}  \end{array} \right].
\end{align*}
Therefore, extracting $\phi_{i,t}$ for each $i\in [n]$ from the preceding relation completes the proof. $\hfill \blacksquare $\\

\noindent
\textbf{\emph{Proof of Lemma \ref{mixture}}}. Since the network is strongly connected and the corresponding $W$ is irreducible and aperiodic, by standard properties of stochastic matrices (see e.g. \cite{rosenthal1995convergence}), the diagonalizable matrix $W$ satisfies
\begin{align}\label{MCP}
\left\|\base_i^\tr W^t-\pi^\tr \right\|_1 \leq n \lambda_{\max}(W)^t, 
\end{align}
for any $i \in [n]$, where $\pi$ is the stationary distribution of a Markov chain with transition kernel $W$. Let us observe the following inequality
\begin{align*}
n \lambda_{\max}(W)^{t-\tau}\leq 2  \ \ \ \ \ \ \ \ \text{for} \ \ \ \ \ \ \ \ t-\tau \geq \tilde{t}\triangleq \frac{\log \left[\frac{n}{2}\right]}{\log \lambda_{\max}(W)^{-1}},
\end{align*}
and recall that the identity $\left\|\base_i^\tr W^{t-\tau}-\pi^\tr\right\|_1\leq 2$ always holds since any power of $W$ is stochastic. With that in mind, we use \eqref{MCP} to break the following sum into two parts to get
\begin{align*}
\sum_{\tau=1}^{t}\sum_{j=1}^n\left|\left[W^{t-\tau}\right]_{ij}-\pi(j)\right| &= \sum_{\tau=1}^{t}\left\|\base_i^\tr W^{t-\tau}-\pi^\tr \right\|_1\\
&=\sum_{\tau=1}^{t-\tilde{t}}\left\|\base_i^\tr W^{t-\tau}-\pi^\tr \right\|_1+\sum_{\tau=t-\tilde{t}+1}^{t}\left\|\base_i^\tr W^{t-\tau}-\pi^\tr \right\|_1\\
& \leq \sum_{\tau=1}^{t-\tilde{t}} n \lambda_{\max}(W)^{t-\tau} + 2 \tilde{t}-2\\
& \leq \frac{n\lambda_{\max}(W)^{\tilde{t}}}{1-\lambda_{\max}(W)} + 2 \tilde{t} \\
&=\frac{2}{1-\lambda_{\max}(W)}+\frac{2\log \left[\frac{n}{2}\right]}{\log \lambda_{\max}(W)^{-1}},
\end{align*}
for any $i\in [n]$. Noting that $1-\lambda_{\max}(W) \leq \log \lambda_{\max}(W)^{-1}$, we have
\begin{align*}
\sum_{\tau=1}^{t}\sum_{j=1}^n\left|\left[W^{t-\tau}\right]_{ij}-\pi(j)\right| \leq \frac{2+2\log \left[\frac{n}{2}\right]}{1-\lambda_{\max}(W)} \leq \frac{4\log n}{1-\lambda_{\max}(W)}.
\end{align*}$\hfill \blacksquare $\\

\noindent
We use the following inequality in \cite{mcdiarmid1998concentration} in the proof of Lemma \ref{Distributed Beliefs}. 
\noindent
\begin{lemma}{\bf (McDiarmid's Inequality)}\label{Mc} 
Let $X_1,...,X_N \in \chi$ be independent random variables and consider the mapping $H: \chi^N \mapsto \Real$. If for $i\in \{1,...,N\}$, and every sample $x_1,...,x_N,x'_i\in \chi$, the function $H$ satisfies
\begin{align*}
\left|H(x_1,...,x_{i-1},x_{i},x_{i+1},...,x_N)-H(x_1,...,x_{i-1},x'_{i},x_{i+1},...,x_N)\right| \leq c_i,
\end{align*}
then for all $\varepsilon>0$,
\begin{align*}
\Prob\bigg\{H(x_1,...,x_N)-\Ex\left[H(X_1,...,X_N)\right] \geq \varepsilon\bigg\}\leq \exp\left\{\frac{-2\varepsilon^2}{\sum_{i=1}^Nc_i^2}\right\}.
\end{align*}
\end{lemma}

\noindent
\textbf{\emph{Proof of Lemma \ref{Distributed Beliefs}}}. According to Lemma \ref{Beliefs}, we have
\begin{align*}
\mu_{i,t}(1)&=\frac{\exp\left\{\eta \phi_{i,t}(1)\right\}}{\sum_{k=1}^m\exp\left\{\eta \phi_{i,t}(k)\right\}}\\
&=\left(1+\sum_{k=2}^m\exp\left\{\eta \phi_{i,t}(k)-\eta \phi_{i,t}(1)\right\}\right)^{-1}\\
&\geq 1-\sum_{k=2}^m\exp\left\{\eta \phi_{i,t}(k)-\eta \phi_{i,t}(1)\right\}, \label{E1}\numberthis
\end{align*}
where we used the fact that $(1+x)^{-1}\geq 1-x$ for any $x \geq 0$. Since we know
\begin{align*}
\|\mu_{i,t}-\base_1\|_{\text{TV}}=\frac{1}{2}\left(1-\mu_{i,t}(1)+\sum_{k=2}^m\mu_{i,t}(k)\right)=1-\mu_{i,t}(1),
\end{align*}
we can simplify \eqref{E1} as follows
\begin{align}
\|\mu_{i,t}-\base_1\|_{\text{TV}} \leq \sum_{k=2}^m\exp\left\{\eta \phi_{i,t}(k)-\eta \phi_{i,t}(1)\right\}.\label{E2}
\end{align} 
For any $k\in [m]$, define
\begin{align*}
\Phi_{i,t}(k)\triangleq \sum_{\tau=1}^{t}\sum_{j=1}^n\left[W^{t-\tau}\right]_{ij}\log \ell_j(\cdot |\theta_k),
\end{align*}
and note that $\Phi_{i,t}(k)$ is a function of $nt$ random variables. As required in McDiarmid's inequality in Lemma \ref{Mc}, set $H=\Phi_{i,t}(k)$, fix the samples for $nt-1$ random variables, and draw two different samples $s_{j,\tau}$ and $s'_{j,\tau}$ for some $j\in [n]$ and some $\tau \in [t]$. The fixed samples are simply cancelled in the subtraction, and we have 
\begin{align*}
|H(...,s_{j,\tau},...)-H(...,s'_{j,\tau},...)|=\left|\left[W^{t-\tau}\right]_{ij}\left(\log \ell_j(s_{j,t} |\theta_k)-\log \ell_j(s'_{j,t} |\theta_k)\right)\right|\leq \left[W^{t-\tau}\right]_{ij} 2B,
\end{align*}
where we used assumption {\bf A1}. Since any power of $W$ is stochastic, summing over $j\in [n]$ and $\tau \in [t]$, we get
\begin{align*}
\sum_{\tau=1}^{t}\sum_{j=1}^n\left(\left[W^{t-\tau}\right]_{ij} 2B\right)^2 \leq 4B^2t.
\end{align*}
We now apply McDiarmid's inequality in Lemma \ref{Mc} to obtain
\begin{align*}
\Prob\big(\phi_{i,t}(k)-\phi_{i,t}(1)>\Ex\left[\Phi_{i,t}(k)\right]-\Ex\left[\Phi_{i,t}(1)\right]+\varepsilon \big) \leq \exp\left\{\frac{-\varepsilon^2}{2B^2t}\right\},
\end{align*}
for $k=2,...,m$. Setting the probability above to $\delta/m$ and taking a union bound over all states, we have for any $k=2,...,m$
\begin{align}
\Prob\left(\phi_{i,t}(k)-\phi_{i,t}(1)\leq \Ex\left[\Phi_{i,t}(k)\right]-\Ex\left[\Phi_{i,t}(1)\right]+\sqrt{2B^2t\log \frac{m}{\delta} }\right) \geq 1-\delta. \label{E7}
\end{align}
On the other hand, in view of assumption {\bf A1}, we have
\begin{align*}
\Ex\left[\Phi_{i,t}(k)-\Phi_{i,t}(1)\right] &=\sum_{\tau=1}^{t}\sum_{j=1}^n\left[W^{t-\tau}\right]_{ij}\Ex\left[\log \ell_j(\cdot |\theta_k)-\log \ell_j(\cdot |\theta_1)\right]\\
&=\sum_{\tau=1}^{t}\sum_{j=1}^n\left(\left[W^{t-\tau}\right]_{ij}-\pi(j)\right)\Ex\left[\log \ell_j(\cdot |\theta_k)-\log \ell_j(\cdot |\theta_1)\right]\\
&~~~~~~~~~~~~~~~~~~~~+\sum_{\tau=1}^{t}\sum_{j=1}^n \pi(j) \Ex\left[\log \ell_j(\cdot |\theta_k)-\log \ell_j(\cdot |\theta_1)\right]\\
&\leq 2B\sum_{\tau=1}^{t}\sum_{j=1}^n\left|\left[W^{t-\tau}\right]_{ij}-\pi(j)\right|-t\sum_{j=1}^n \pi(j)D_{KL}\left(\ell_j(\cdot |\theta_1)\| \ell_j(\cdot |\theta_k)\right)\\
&= 2B\sum_{\tau=1}^{t}\sum_{j=1}^n\left|\left[W^{t-\tau}\right]_{ij}-\pi(j)\right|-\mathcal{I}(\theta_1,\theta_k)t\\
&\leq \frac{8B\log n}{1-\lambda_{\max}(W)}-\mathcal{I}(\theta_1,\theta_k)t,
\end{align*}
where we applied Lemma \ref{mixture} to derive the last step. Using \eqref{orderassum}, we simplify above to get 
\begin{align}
\Ex\left[\Phi_{i,t}(k)-\Phi_{i,t}(1)\right]&\leq \frac{8B\log n}{1-\lambda_{\max}(W)}-\mathcal{I}(\theta_1,\theta_2)t, \label{E5}
\end{align}
for any $k=2,...,m$. Plugging \eqref{E5} into \eqref{E7} and combining with \eqref{E2}, we have {\small
\begin{align*}
\|\mu_{i,t}-\base_1\|_{\text{TV}} &\leq \sum_{k=2}^m \exp \left\{-\eta\mathcal{I}(\theta_1,\theta_2)t+\eta \sqrt{2B^2t\log \frac{m}{\delta} }+\frac{8\eta B\log n}{1-\lambda_{\max}(W)}\right\}\\
&\leq m \exp \left\{-\eta\mathcal{I}(\theta_1,\theta_2)t+\eta \sqrt{2B^2t\log \frac{m}{\delta} }+\frac{8\eta B\log n}{1-\lambda_{\max}(W)}\right\},
\end{align*} }
with probability at least $1-\delta$, and thereby completing the proof. $\hfill \blacksquare $\\

\noindent
\textbf{\emph{Proof of Lemma \ref{DecRegretDecay}}}. 
We recall from the statement of the lemma that $q_{i,t}(k)=\phi_{i,t}(k)-\phi_{t}(k)$, and calculate the ratio $\mu_{i,t}(k)/\mu_{t}(k)$ for any $k\in [m]$ as follows,
\begin{align*}
\frac{\mu_{i,t}(k)}{\mu_t(k)}&=\exp\left\{\eta q_{i,t}(k)\right\}\frac{\Ex_{\mu_0}\left[\exp\left\{\eta \phi_{t}\right\}\right]}{\Ex_{\mu_0}\left[\exp\left\{\eta \phi_{i,t}\right\}\right]}\\
&=\exp\left\{\eta q_{i,t}(k)\right\}\frac{\Ex_{\mu_0}\left[\exp\left\{\eta \phi_{t}\right\}\right]}{\Ex_{\mu_0}\left[\exp\left\{\eta \phi_{t}\right\}\exp\left\{\eta q_{i,t}\right\}\right]}\\
&=\exp\left\{\eta q_{i,t}(k)\right\}\frac{1}{\Ex_{\mu_0}\left[\frac{\exp\left\{\eta \phi_{t}\right\}}{\Ex_{\mu_0}\left[\exp\left\{\eta \phi_{t}\right\}\right]}\exp\left\{\eta q_{i,t}\right\}\right]}\\
&=\exp\left\{\eta q_{i,t}(k)\right\}\frac{1}{\Ex_{\mu_0}\left[\frac{\mu_{t}}{\mu_0}\exp\left\{\eta q_{i,t}\right\}\right]}\\
&=\exp\left\{\eta q_{i,t}(k)\right\}\frac{1}{\Ex_{\mu_{t}}\left[\exp\left\{\eta q_{i,t}\right\}\right]}.
\end{align*}
This entails
\begin{align*}
\frac{1}{\eta}\Ex_{\mu_{i,t}}\left[\log\frac{\mu_{i,t}}{\mu_{t}}\right]&=\Ex_{\mu_{i,t}}\left[q_{i,t}\right]-\frac{1}{\eta}\log \Ex_{\mu_{t}}\left[\exp\left\{\eta q_{i,t}\right\}\right] \leq \Ex_{\mu_{i,t}}\left[q_{i,t}\right]-\Ex_{\mu_{t}}\left[q_{i,t}\right],
\end{align*}
where we used Jensen's inequality on the convex function $-\log(\cdot)$. Setting the expectation measures in the right hand side of above to $\mu_t$, and recalling the ratio $\mu_{i,t}/\mu_{t}$ from above, we conclude that,
\begin{align*}
\frac{1}{\eta}\Ex_{\mu_{i,t}}\left[\log\frac{\mu_{i,t}}{\mu_{t}}\right]&\leq \Ex_{\mu_t}\left[\frac{\mu_{i,t}}{\mu_{t}}q_{i,t}\right]-\Ex_{\mu_t}\left[q_{i,t}\right]\\
&= \Ex_{\mu_t}\left[\left(\frac{\exp\{\eta q_{i,t}\}}{\Ex_{\mu_{t}}\left[\exp\{\eta q_{i,t}\} \right]}-1\right)q_{i,t}\right]\\
&= \Ex_{\mu_t}\left[\bigg(\frac{\exp\{\eta q_{i,t}\}}{\Ex_{\mu_{t}}\left[\exp\{\eta q_{i,t}\} \right]}-1\bigg)\bigg(q_{i,t}-\Ex_{\mu_t}[q_{i,t}]\bigg)\right]\\
&\leq \sqrt{\Ex_{\mu_t}\left[\left(\frac{\exp\{\eta q_{i,t}\}}{\Ex_{\mu_{t}}\left[\exp\{\eta q_{i,t}\} \right]}-1\right)^2\right]\bigg(\V_{\mu_t} [q_{i,t}]\bigg)}, \numberthis \label{RBound C5}
\end{align*}
where we applied Cauchy-Schwarz inequality in the last line. Then, we appeal to Jensen's inequality again to get
\begin{align*}
\Ex_{\mu_t}\left[\left(\frac{\exp\{\eta q_{i,t}\}}{\Ex_{\mu_{t}}\left[\exp\{\eta q_{i,t}\} \right]}-1\right)^2\right]&=\Ex_{\mu_t}\left[\left(\frac{\exp\{\eta q_{i,t}\}}{\Ex_{\mu_{t}}\left[\exp\{\eta q_{i,t}\} \right]}\right)^2\right]-1\\
&\leq\Ex_{\mu_t}\left[\left(\frac{\exp\{\eta q_{i,t}\}}{\exp\{\Ex_{\mu_{t}}\left[\eta q_{i,t}\right]\} }\right)^2\right]-1\\
&=\Ex_{\mu_t}\bigg[\exp\bigg\{2\eta\bigg(q_{i,t}-\Ex_{\mu_{t}}\left[q_{i,t}\right]\bigg)\bigg\}\bigg]-1.
\end{align*}
Note that the function $g(z)=(\exp\{z\}-1-z)/z^2$ is nondecreasing over reals, and let $z=2\eta(q_{i,t}-\Ex_{\mu_{t}}\left[q_{i,t}\right])$ in $g(z)$. The condition $\eta\|q_{i,t}\|_{\infty} \leq 1/4$ immediately implies that $z \leq 1$, so recalling that $z$ is the argument of exponential in above, we bound the right hand side as,
\begin{align*}
\Ex_{\mu_t}\bigg[\exp\bigg\{2\eta\bigg(q_{i,t}-\Ex_{\mu_{t}}\left[q_{i,t}\right]\bigg)\bigg\}\bigg]-1\leq 4\left(\exp(1)-2\right)\V_{\mu_t}\left[\eta q_{i,t}\right] \leq 4\V_{\mu_t}\left[\eta q_{i,t}\right].
\end{align*}
Plugging the bound above into \eqref{RBound C5}, results in
\begin{align*}
\frac{1}{\eta}\Ex_{\mu_{i,t}}\left[\log\frac{\mu_{i,t}}{\mu_{t}}\right]&\leq \sqrt{4\V_{\mu_t}\left[q_{i,t}\right]\V_{\mu_t}\left[\eta q_{i,t}\right]}= 2\eta \V_{\mu_t}\left[q_{i,t}\right].
\end{align*}
Summing above over $t\in [T]$ and recalling \eqref{Regret}, concludes the proof.$\hfill \blacksquare $\\


\bibliographystyle{IEEEtran}
\bibliography{IEEEabrv,shahin}

\end{document}